\documentstyle[11pt]{article}
 
\newtheorem{thm}{Theorem}[section]
\newtheorem{prop}[thm]{Proposition}

\newtheorem{lemma}[thm]{Lemma}

\newcounter{ex}[section]

\newcommand{\red}{{\rm red}}

\newcommand{\E}{{\cal E}}

\newcommand{\C}{{\bf C}}

\newcommand{\lo}{{\ \rightarrow\ }}

\newcommand{\R}{{\bf R}}
\newcommand{\Q}{{\bf Q}}

\newcommand{\ep}{\epsilon}

\newcommand{\Z}{{\bf Z}}
\renewcommand\d{\delta}
\newcommand{\F}{{\cal F}}
\newcommand{\ti}{\tilde}
\newcommand{\Spec}{{\rm Spec}}
\renewcommand{\O}{{\cal O}}

\newcommand{\omlog}{\Omega^1_{X/\Z}(\log X^{\red}_S/\log S)}

\newcommand{\lon}{{\longrightarrow}}

\newenvironment{eq}{\addtocounter{thm}{1}\begin{equation}
}{\end{equation}}
\newenvironment{eqnar}{\addtocounter{thm}{1}\begin{eqnarray}
}{\end{eqnarray}}

\newcommand{\Om}{\Omega}

\newcommand{\mod}{{\rm mod}}

\newcommand{\noi}{\noindent}

\newcommand{\hatch}{{\widehat {\rm ch}}}
\newcommand{\hatTd}{{\widehat {\rm Td}}}
\newcommand{\hCH}{{\widehat {\rm CH}}}
\newcommand{\CH}{{\rm CH}}
\newcommand{\ch}{{\rm ch}}
\newcommand{\tch}{\widetilde {\rm ch}}
\newcommand{\tTd}{\widetilde {\rm Td}}
\newcommand{\Gr}{{\rm Gr}}
\newcommand{\hch}{{\widehat {\rm ch}}}
\newcommand{\Td}{{\rm Td}}
\newcommand{\hK}{\widehat {\rm K}}
\newcommand{\pione}{{\bf P}^1}

\begin{document}
\title{$\ep$-constants and Arakelov Euler characteristics}

\author{Ted Chinburg\thanks{Supported in part by NSF grant DMS97-01411.}, Georgios Pappas\thanks{Supported in part by NSF grant DMS99-70378 and by a Sloan Research Fellowship.} and Martin J. Taylor\thanks{EPSRC Senior Research Fellow.}}\date{}
\maketitle
 \section{Introduction} \label{intro}

Let $X$ be a regular scheme projective and flat over $\Spec(\Z)$,
equidimensional of relative dimension $d$. Consider the
Hasse-Weil zeta function of $X$, $\zeta(X, s)=\prod_x(1-N(x)^{-s})^{-1}$
where $x$ ranges over the closed points of $X$ and $N(x)$ is the order of the residue field
of $x$. Denote by $L(X, s)$ the zeta function with $\Gamma$-factors
$L(X, s)=\zeta(X, s)\Gamma(X, s)$. The $L$-function conjecturally satisfies
a functional equation
$$
L(X, s)=\epsilon(X)A(X)^{-s}L(X, d+1-s)
$$
where $\ep(X)$ and $A(X)$ are real numbers defined independently
of any conjectures (the ``$\ep$-constant" and the ``conductor").
In fact, the unconditional definition of $\ep(X)$ and $A(X)$
involve choices of auxiliary primes $l$ with embeddings
$\Q_l\subset\C$ (see [De]). In this note, we will suppress
any notation regarding these choices; this should not cause
any confusion.
 
The purpose of this note is to explain a way to obtain the absolute value 
$|\ep(X)|$ as an ``arithmetic" Euler de Rham characteristic in the framework 
of the higher dimensional Arakelov theory of Gillet and Soul\'e.
Choose a hermitian metric on the tangent bundle of $X(\C)$
which is K\"ahler; it gives a hermitian metric on $\Omega^1_{X_\C}$.
Recall the definition of the arithmetic Grothendieck group $\hK_0(X)$
of hermitian vector bundles
of Gillet and Soul\'e ([GS1, II, \S 6]; all
hermitian metrics are smooth and invariant under the complex conjugation
on $X(\C)$).
There is an arithmetic Euler characteristic homomorphism
$$
\chi_Q: \hK_0(X)\ \lon \ \R
$$
such that if $(\F, h)$ is a vector bundle on $X$ with a hermitian metric
on $\F_\C$, then $\chi_Q((\F, h))$ is the Arakelov degree of the
hermitian line bundle on $\Spec(\Z)$ formed by the determinant of the cohomology
of $\F$ with its Quillen metric.
The arithmetic Grothendieck group $\hK_0(X)$ is a $\lambda$-ring
with $\lambda^i$-operations defined in loc. cit. \S 7:
If $(\F, h)$ is the class of a vector bundle with a hermitian metric on $\F_\C$
then $\lambda^i((\F, h))$ is the class of the vector bundle $\wedge^i\F$
with the exterior power metric on $\wedge^i\F_\C$ induced from $h$. 
Now consider the sheaf of differentials $\Omega^1_{X/\Z}$;
this is a ``hermitian coherent sheaf" in the terminology
of [GS3, 2.5]. Since $X$ is regular, by loc. cit. 2.5.2
$\Omega^1_{X/\Z}$ defines an element $\Omega$ in $\hK_0(X)$ as follows:
Each embedding of $X$ into projective space over 
$\Spec(\Z)$ gives a short exact
sequence
$$
\E\ :\ 0\ \lo\ N\ \lo\ P\ \lo\ \Omega^1_{X/\Z}\ \lo\ 0
$$
with $P$ and $N$ vector bundles on $X$
(here $P$ is the restriction of the relative differentials
of the projective space to $X$ and $N$ is the conormal bundle of 
the embedding). Pick hermitian metrics 
$h^P$ and $h^N$ on $P_\C$ and $N_\C$ respectively
and denote by  ${\tch}(\E_\C)$ the secondary
Bott-Chern characteristic class of the exact
sequence of hermitian vector bundles $\E_\C$
(as defined in [GS1]; there is a difference 
of a sign between
this definition and the definition in [GS3, 2.5.2]).
Then
$$
\Omega=((P, h^P), 0)-((N, h^N), 0)+((0,0), {\tch}(\E_\C))\in \hK_0(X)
$$ 
depends only on the original choice of
K\"ahler metric.  

For each $i\geq 0$ we can
consider now the element
$\lambda^i(\Omega)$ in $\hK_0(X)$.
Motivated by the ``higher dimensional Fr\"ohlich
conjecture" of [CEPT], we conjecture that 
\begin{eq} \label{con}
-\log|\ep(X)|=\sum_{i=0}^d(-1)^i{\chi_Q}(\lambda^i(\Omega)).
\end{eq}
Denote by $X_S$ the disjoint union of the 
singular fibers of $f:X\lo \Spec(\Z)$.
In [B], S. Bloch conjectures that the conductor $A(X)$
is given by
$$
A(X)={\rm ord}((-1)^dc^{X_S}_{d+1}(\Omega^1_{X/\Z}))
$$
where $c_{d+1}^{X_S}(\Omega^1_{X/\Z}):=c_{d+1, X_S}^{X}(\Omega^1_{X/\Z})\cap [X]$ is the
localized $d+1$-st Chern class in $\CH_0(X_S)$
described in loc. cit. Here for a zero cycle
$\sum_in_ix_i$, ${\rm ord}(\sum_in_ix_i)=\prod_i(\# k(x_i))^{n_i}$,
with $k(x_i)$ the residue field of $x_i$. In this paper we show:

\begin{thm} \label{main} 
The equality \ref{con} is equivalent to Bloch's conjecture.
\end{thm}

The main ingredients in the proof 
are the Arithmetic Riemann-Roch theorem
of Gillet and Soul\'e and the fact (Proposition \ref{locAr})
that Bloch's localized Chern class agrees
with the corresponding ``arithmetic" 
Chern class of Gillet-Soul\'e.

Since Bloch has proven in [B] his conjecture
for an arithmetic surface ($d=1$) we see
that \ref{con} holds in this case.
In this note we also show:

\begin{thm} \label{main2}
Bloch's conjecture, and 
therefore equality \ref{con}, holds when for all primes $p$, the fiber of $X \lo \Spec(\Z)$ over $p$
is a divisor with strict normal crossings with multiplicities relative
prime to $p$.
\end{thm}

In fact, under the hypothesis of the above theorem,
we can show \ref{con} directly by replacing the use of the
arithmetic Riemann-Roch theorem by Serre duality
and the fact, due to Ray and Singer ([RS], Theorem 3.1),
that the analytic torsion of the de Rham complex 
is trivial. We are grateful to C. Soul\'e for 
pointing this out to us; this approach is explained
in detail in [CPT2]. Also, since we have $\ep(X)^2=A(X)^{d+1}$, 
we could have expressed
\ref{con} using the conductor $A(X)$. However, it seems that \ref{con}
is more canonical and it could generalize in a motivic framework (for
example to varieties with a group action).
Indeed, the inspiration for \ref{con} comes from [CEPT], see also [CPT1], where
we observed a close connection between an equivariant version of an Euler de 
Rham characteristic as above and $\ep$-constants. Viewed
this way, Theorem \ref{main} also provides some  
indirect positive evidence for 
the general higher dimensional Fr\"ohlich conjecture of [CEPT]. 
In [CPT2], we use the results of this note 
to obtain the actual $\ep$-constant (not just its absolute value) 
of the Artin motive obtained from the pair $(X, V)$ of an arithmetic variety 
$X$ with an action of a finite group $G$ and a 
symplectic character $V$ of $G$. 

We would like to express our thanks to C. Soul\'e;
this note would not have existed without his advice.
We would also
like to thank T. Saito for useful conversations and B. Erez for pointing out
the reference [A].
After a preliminary version of this note was completed we have
learned that K. Kato and T. Saito have announced a proof
of a stronger version of Theorem \ref{main2} in which the
assumption on the multiplicities is dropped; their proof
is significantly more involved than the proof
of the tame case that we consider here.
T. Saito informed us that a similar argument to ours for the proof
of the tame case is given by K. Arai in his thesis,
which is currently in preparation.

\section{Arithmetic Riemann-Roch} \label{calARR}

The formulae of [De] imply that $\ep(X)^{2}=A(X)^{d+1}$ 
(we can see that this also follows directly from the 
conjectural functional equation).
Therefore, \ref{con}
translates to
\begin{eq} \label{con2}
{{d+1}\over 2}\cdot \log A(X)= -\sum_{i=0}^d(-1)^i{\chi_Q}(\lambda^i(\Omega)).
\end{eq}
Denote by $\hCH^{\cdot}(X)$, $\hCH_{\cdot}(X)$ the arithmetic Chow groups
of Gillet and Soul\'e ([GS1-2]), graded by codimension and 
dimension of cycles respectively. Since $X_S$ has empty generic fiber, 
there is a natural homomorphism
$$
z_S: \CH_0(X_S)\lo \hCH_0(X)=\hCH^{d+1}(X).
$$
The direct image homomorphism 
$$
f_*: \hCH^{d+1}(X)\lo \hCH^1(\Spec(\Z))=\R
$$
satisfies
$
f_*(z_S(a))=\log({\rm ord}(a))$ for $a\in \CH_0(X_S)$.
Therefore, Theorem \ref{main} will follow if we show:
\begin{thm} \label{ARR} 
$\ \ \ \displaystyle{\sum_{i=0}^d(-1)^i{\chi_Q}(\lambda^i(\Omega))=(-1)^{d+1}{d+1\over 2}f_*(z_S(c^{X_S}_{d+1}(\Omega^1_{X/\Z}))}).$
\end{thm}

In what follows we will use heavily the notations 
and results of
[GS1], [GS2] and [GS3].

First observe that from the definition of $\Omega$,
we obtain $\ch(\Omega)=\ch(\Omega^1_{X(\C)})$,
where $\ch$ denotes the Chern character form
(its domain can be extended to $\hK_0(X)$
as in [GS1]). By [GS1, Lemma 7.3.3], we
have $\ch(\Omega^i_{X(\C)})=\lambda^i(\ch(\Omega^1_{X(\C)}))$;
here $\Omega^i_{X(\C)}$, $0\leq i\leq d$, 
has the exterior power metric and the
$\lambda$-ring structure on differential forms
is given by the grading as in loc. cit.  
We obtain that $\ch(\lambda^i(\Omega))
=\lambda^i(\ch(\Omega))=\lambda^i(\ch(\Omega^1_{X(\C)}))=\ch(\Omega^i_{X(\C)})$
where the first equality follows
from the fact that $\ch=\omega\cdot\hch:
\hK_0(X)\lo A(X_\R)$ is a $\lambda$-ring homomorphism
(see loc. cit.).
 
From the Arithmetic Riemann Roch theorem of Gillet and Soul\'e ([GS3], Theorem 7,
see also 4.1.5 loc. cit.) we 
now have
\begin{eqnarray} 
\sum_{i=0}^d(-1)^i{\chi_Q}(\lambda^i(\Omega))&=&
f_*((\hatch (\sum_{i=0}^d(-1)^i\lambda^i(\Omega))\cdot {\hatTd}(X))^{(d+1)})\label{ARReq}\\
&&-{1\over 2}\int_{X(\C)}\ch(\sum_{i=0}^d(-1)^i\Om^i_{X(\C)}){\rm Td}(T_{X(\C)})R(T_{X(\C)})\nonumber
\end{eqnarray}
where the notations are as in loc. cit.
and the factor of $1/2$
in front of the second term results
from the normalization discussed after
equation (15) in section 4.1.5. We
first show:
\begin{prop} \label{int}
$\ \ \displaystyle{
\int_{X(\C)}\ch(\sum_{i=0}^d(-1)^i\Om^i_{X(\C)}){\rm Td}(T_{X(\C)})R(T_{X(\C)})=0}.
$
\end{prop}

\smallskip

\begin{Proof} (Shown to us by C. Soul\'e.) By the classical identity
applied on the level of Chern forms we obtain
$$
\ch(\lambda_{-1}(\Om^1_{X(\C)}))\Td(T_{X(\C)}))=c_{d}(T_{X(\C)})
$$
(see [R, 6.19]). Therefore
the integral is equal to:
$$
\int_{X(\C)}c_d(T_{X(\C)})R(T_{X(\C)}).
$$
But $R(T_{X(\C)})$ is non-zero in positive degrees only;
therefore the degree of the form $c_d(T_{X(\C)})R(T_{X(\C)})$
is at least $d+1$ and the integral vanishes.
\end{Proof}

\medskip
It remains to deal with the first term of the 
right hand side of \ref{ARReq}. We will show:
\begin{prop} \label{chTd} 
$$(\hatch (\sum_{i=0}^d(-1)^i\lambda^i(\Omega))\cdot {\hatTd}(X))^{(d+1)}
=(-1)^{d+1}{d+1\over 2}\hat c_{d+1}(\Omega).$$
\end{prop}

\smallskip

\begin{Proof}
Recall the definition of $\hatTd(X)$
from [GS3];
we have an exact sequence
$$
\E^*_\C: 0\lo T_{X_\C}=({\Omega^{1}_{X_\C}})^*\lo P^*_\C\lo N^*_\C\lo 0.
$$
We set
$$
\hatTd(X):=\hatTd(\bar P^*)\hatTd^{-1}(\bar N^*)+a(\tTd(\E^*_\C)\Td(\bar N^*_\C)^{-1})
$$
where $\tTd(\E^*_\C)$ is the Todd-Bott-Chern secondary 
form attached to the sequence $\E^*_\C$ (see [GS3], p. 503)
and $\Td$ is the usual Todd form. We are just
interested in the terms of degree $0$ and $1$ of $\hatTd(X)$.
If $\bar E$ is a hermitian vector bundle, we have
$$
\hatTd(\bar E^*)=1+{\hat c_1(\bar E^*)\over 2}+\cdots,\quad
\hatTd^{-1}(\bar E^*)=1-{\hat c_1(\bar E^*)\over 2}+\cdots
$$ 
The $(0,0)$ component of $\Td(\bar N^*_\C)^{-1}$ is $1$.
We can also see that the $(0,0)$ component of the secondary form
$\tTd(\E^*_\C)$ is given by
$$
\tTd(\E^*_\C)^{(0,0)}={\ti c_1(\E^*_\C)\over 2}
$$
where $\ti c_1(\E^*_\C)$
the ``secondary" first Bott-Chern 
form associated to $\E^*_\C$.
This gives
$$
\hatTd(X)=1+{\hat c_1(\bar P^*)-\hat c_1(\bar N^*)\over 2}+a({\ti c_1(\E^*_\C)\over 2})+\cdots=
1+{\hat c_1(\Omega^*)\over 2}+\cdots 
$$
and therefore
\begin{eq} \label{congr}
\hatTd(X)=\hatTd(\Omega^*)\ \mod\ \hCH^{\geq 2}(X)_\Q.
\end{eq}

Now consider the $\gamma$ operations on the $\lambda$-ring
$\hK_0(X)$ with augmentation $\ep: \hK_0(X)\lo \Z$ 
given by $\ep((\bar E, \eta))={\rm rk}(E)$ (see [R, \S 4]).
If $\ep(x)=d$, then (as in [CPT1] \S 1) we have:
\begin{eq}\label{gala}
(-1)^d\gamma^d(x-\ep(x))=\sum_{i=0}^d(-1)^i\lambda^i(x).
\end{eq}

Therefore $\hatch(\sum_{i=0}^d(-1)^i\lambda^i(x))$
is concentrated in degrees $d$ and $d+1$ only
and so by \ref{congr}
$$
\hatch (\sum_{i=0}^d(-1)^i\lambda^i(\Omega))\cdot {\hatTd}(X)=
\hatch (\sum_{i=0}^d(-1)^i\lambda^i(\Omega))\cdot {\hatTd}(\Omega^*).
$$

By the above and \ref{gala}
it is enough to show that for $x\in \hK_0(X)$ we have
$$
(\hatch (\gamma^d(x-\ep(x))\cdot {\hatTd}(x^*))^{(d+1)}=-{d+1\over 2}
\hat c_{d+1}(x).
$$

Let $a_1, \ldots , a_{d+1}$ be
the ``arithmetic Chern roots" of $x$. By definition, these
are formal symbols
such that the 
arithmetic Chern classes of $x$ are the elementary 
symmetric functions of $a_i$; we can perform
our calculation using these symbols. The Chern roots of the dual $x^*$
are $-a_1,\ldots , -a_{d+1}$. A standard argument using [GS1, Theorem 4.1]
shows that we have 
$$
\hatch(\gamma^d(x-\ep(x)))=\sum^{d+1}_{i=0}\prod_{j\neq i}(e^{a_j}-1),
$$
while by definition
$$
\hatTd (x^*)=\prod_{i=0}^{d+1}{{-a_i}\over {1-e^{-(-a_i)}}}=\prod_{i=0}^{d+1}{{a_i}\over {e^{a_i}-1}}.
$$

The product is equal to
$$
\sum_{j=1}^{d+1}{a_1a_2\cdots a_{d+1}\over e^{a_j}-1}=\sum_{j=1}^{d+1}
(a_1\cdots \widehat a_j\cdots a_{d+1}-{a_1\cdots a_{d+1}\over 2})+\cdots
$$
$$
=\hat c_d(x)-{d+1\over 2}\hat c_{d+1}(x)+\cdots
$$
which gives the desired result.
\end{Proof}

\medskip
\medskip
\section{Localized Chern classes.} \label{locCh}

We continue with the same assumptions and notations.
Recall the homomorphism
$$
z_S: \CH_0(X_S)_\Q\lo \hCH_0(X)_\Q=\hCH^{d+1}(X)_\Q.$$

\begin{prop} \label{locAr}
$\ \
z_S(c^{X_S}_{d+1}(\Omega^1_{X/\Z}))=\hat c_{d+1}(\Omega).
$
\end{prop}

\medskip

Theorem \ref{ARR} follows from Propositions
\ref{locAr}, \ref{int}, \ref{chTd} and equation \ref{ARReq}.

\medskip

\noi{\bf Proof of Proposition \ref{locAr}:}
We review the construction of
the localized Chern class via the Grassmannian graph construction (as described in [B]
\S 1, or in [GS3] \S 1) applied to the complex
$0\lo N\buildrel \d\over \lo P$ with cokernel $\Omega^1_{X/\Z}$.
Set $U=X-X_S$.
Let $p$ be the projection $X\times \pione \lo X$.
Set $M:=p^*N(1)\oplus p^*P$
where $(1)$ denotes the Serre twist (which we view as tensoring
with the pull-back of $\O_{\pione}(\infty)$ under $X\times\pione\lo \pione$). Let us consider 
the Grassmannian $\Gr(r, M)$ over $X\times \pione$ of rank $r={\rm rk}(N)$ 
local direct summands 
of $M$. Denote by $\pi_0: \Gr(r, M)\lo X\times \pione$ the natural projection
morphism. The diagonal embedding $p^*N\subset p^*N(1)\oplus p^*P$ gives
a section $s$ of $\pi_0$ over the subscheme $(X\times{\bf A}^1)\cup
(U\times \pione)$. In fact, over $X\times {\bf A}^1$ the image
of $p^*N$ can be identified with the graph of $\d$.
Denote by $W$ the Zariski closure of the image 
$$
s((X\times{\bf A}^1)\cup
(U\times \pione))\subset \Gr(r, M);
$$ 
this is an integral
subscheme of $\Gr(r, M)$ which is called the Grassmannian graph of $N \lo P$.
The morphism $\pi:={\pi_0}_{|W}$ is projective and
gives an isomorphism on the generic fibers.
Let $W_{\infty}$ be the effective Cartier divisor on $W$
given by the inverse image of $X\times\{\infty\}$ under $\pi$.
Also let $\tilde X$ be the Zariski closure in $W_\infty$ 
of the restriction of the section $s$ to $U\times \{\infty\}$.
Then $\pi_{|\tilde X}: \tilde X\lo X$ is birational (an isomorphism over
$U$). As in [GS3], we see that the cycle
$$
Z=[W_\infty]-[\tilde X]
$$ 
is supported in the inverse image 
of $X_S$. Looking at supports, we have
$|W_\infty|=|\tilde X|\cup |Z|$.

Denote by
$\xi_1$ the universal subbundle of rank $r$ on $\Gr(r, M)$ and
by $\xi_0$ the ``constant" bundle which is the base change of
$P$ under the (smooth) morphism $\Gr(r, M)\lo X$.  The section $s$
gives
$$
s_\C:X_\C\times \pione_\C=W_\C\lo \Gr(r, M)_\C.
$$
The pull-back of $\xi_0$ under $s_\C$ is $p^*P_\C$; 
the pull-back of $\xi_1$ under $s_\C$ is $p^*N_\C$. 
Denote the restrictions ${\xi_0}_{|W}$, ${\xi_1}_{|W}$
by $\zeta_0$, $\zeta_1$.
Equip ${\zeta_0}_{\C}$, ${\zeta_1}_{\C}$ with the hermitian metrics which correspond
to the hermitian metrics on $p^*P_\C$, $p^*N_\C$ obtained via base change
from the metrics on $P_\C$, $N_\C$. We will denote by $\bar\zeta_1$,
$\bar\zeta_0$ the vector bundles $\zeta_1$, $\zeta_0$
on $W$ endowed with the above hermitian metrics on $W_\C$.
Set $\bar\zeta=(\bar\zeta_0,0)-(\bar\zeta_1,0)\in \hK_0(W)$.

There is a natural morphism $\xi_1\lo \xi_0$ 
obtained by the natural inclusion
$\xi_1\subset \pi^*_0 M$ followed by the projection $\pi^*_0M\lo \xi_0=\pi^*_0p^*P$.
After restricting to $W_\C$ this corresponds 
 to the composition 
$p^*N_\C\lo p^*P_\C$.

Over $X_\C$ we have
the exact sequence
$$
\E_\C \ :\ 0\lo N_\C\lo P_\C\lo \Omega^1_{X_\C}\lo 0.
$$

This gives an exact sequence over $X_\C\times \pione_\C=W_\C$:
$$
p^*\E_\C: 0\lo p^*N_\C\lo p^*P_\C\lo p^*\Omega^1_{X_\C}\lo 0.
$$

Consider 
$A=(0, \tch(p^*\E_\C))$ in $\hK_0(W)$. 
Let us now define the elements
$$
b=\hat c_{d+1}(\bar\zeta+A)\in {\hCH}^{d+1}(W)_\Q,
$$
$$
\mu=\pi_*(b)\in {\hCH}^{d+1}(X\times \pione)_\Q.
$$

\begin{lemma} \label{restr}
The restrictions of $\mu$ to 
$X\times \{0\}$ and $X\times \{\infty\}$ are equal.
\end{lemma}
\begin{Proof}  By [GS2, 
Theorem 4.4.6] the restrictions are well
defined and their difference is given by 
$$
a(\int_{\pione(\C)}\omega(\mu)\log |z|^2)
$$
where $\omega$ and $a$ are defined in [GS2, 3.3.4];
$\omega(\mu)$ is a $(d+1,d+1)$-form on $(X\times\pione)(\C)$
and the integral in the parenthesis gives a $(d,d)$-form on $X(\C)$. 
Since $\pi$ is an isomorphism on the generic fibers,
by the definition of $\bar\zeta$ and $A$, we can see that
the form $\omega(\mu)$ is obtained by 
pulling back via the projection $p_\C:X(\C)\times {\bf P}^1(\C)\lo X(\C)$ 
a $(d+1,d+1)$-form on $X(\C)$.
It follows that
$$
\int_{\pione(\C)}\omega(\mu)\log |z|^2=0
$$
(the integral changes sign when $z$ is replaced by $1/z$).
\end{Proof}
\medskip

Recall that the morphism $\pi :W\lo X\times \pione$ restricts to give
a projective morphism $\pi^{|Z|}: |Z|\lo X_S\times\infty=X_S$. Here
$|Z|$ is the (reduced) support of $Z$.
Set $\xi=\xi_0-\xi_1\in {K}_0(\Gr(r, M))$
and denote by $[Z]$ the fundamental cycle of $Z$
in $\CH_{d+1}(|Z|)$.

\begin{lemma} \label{twores}
a) The restriction of $\mu$ to $X\times\{0\}$ is equal to $\hat c_{d+1}(\Omega)$;

b) The restriction of the class $\mu$ to $X\times\{\infty\}$ is equal to
the image of
$\pi^{|Z|}_*(c_{d+1}(\xi_{||Z|})\cap [Z])\in \CH_0(X_S)_\Q$ under $z_S$. 
\end{lemma}

Before we continue with the proof, let us point out 
that since by definition
$
c_{d+1}^{X_S}(\Omega^1_{X/\Z})=c^{X}_{d+1,X_S}(\Omega^1_{X/\Z})\cap [X]=\pi^{|Z|}_*(c_{d+1}(\xi_{||Z|})\cap [Z]),
$ Lemmas \ref{restr} and \ref{twores} together
imply the proof of Proposition \ref{locAr}.
\medskip 

\begin{Proof}  Part (a) is straightforward; 
indeed $\bar\xi_0$ restricts to give $\bar P$, $\bar \xi_1$
gives $\bar N$ and $A$ gives $(0, \tch(\E_\C))$. 

Let us show part (b). Recall $W$ is integral of dimension $d+2$, $W_\infty$ 
is an effective Cartier divisor in $W$ and we have
$[W_\infty]=Z+[\tilde X]$. Denote by $|W_\infty|$ the reduced support of 
$W_\infty$ in $W$. Since $\pi^{|W_\infty|}: |W_\infty|\lo X\times\{\infty\}=X$
is a projective morphism which is an isomorphism on the generic fiber,
$$
\pi^{|W_\infty|}_*: \hCH^{d+1}(|W_\infty|)_\Q\lo \hCH^{d+1}(X)_\Q
$$
is well-defined. Also, since $i: W_\infty\lo W$ is the inclusion
of an effective Cartier divisor with smooth generic fiber, the pull-back 
$i^*(b)$ makes sense in $\hCH^{d+1}(W_{\infty})_\Q=\hCH^{d+1}(|W_\infty|)_\Q=\hCH_0(|W_\infty|)_\Q$ 
and we have
$$
\mu_{|X\times\{\infty\}}=\pi_*(b)_{|X\times\{\infty\}}=\pi^{|W_\infty|}_*(i^*(b)).
$$
(see for example [GS3, 2.2.7]). 

In what follows, we will calculate $i^*(b)$.
For simplicity set $G=\Gr(r, M)$. Equip the
bundles $\xi_1$, $\xi_0$ on $G$ with hermitian metrics
and set $\bar\xi=(\bar\xi_0, 0)-(\bar\xi_1, 0)\in \hK_0(G)$.
Consider $B=\hat c_{d+1}(\bar\xi)$
in $\hCH^{d+1}(G)$ and 
$B_{|W}=\hat c_{d+1}(\bar\xi_{|W})$ in $\hCH^{d+1}(W)$.
Note that  $\bar\xi_{|W}\in \hK_0(W)$ need not agree
with $\bar\zeta$ because the metrics might not agree. 
In any case, we
can write
\begin{eq}\label{beforecrap}
B_{|W}-\hat c_{d+1}(\bar\zeta+A)=a(\eta)
\end{eq}
with $\eta$ a $(d,d)$-form on $W(\C)$. The pull-back $i^*(B_{|W})$ 
is the $d+1$-st arithmetic Chern class of the
restriction of the bundle $\bar\xi_0-\bar\xi_1$ to $W_\infty$.
We have 
\begin{eq}\label{crap}
i^*(b)=i^*(B_{|W})-a(i^*_\C(\eta)).
\end{eq}

By [GS3, Theorem 4 (1)],
$i^*(B_{|W})=B\cdot_j[W_\infty]$ in $\hCH^{d+1}(|W_\infty|)_\Q=\hCH_0(|W_\infty|)_\Q$;
here $j: |W_\infty|\lo G$ is the natural embedding
and $[W_\infty]\in \hCH_{d+1}(|W_\infty|)_\Q$ is
the fundamental cycle of $W_\infty$ (the notations
are as in loc.cit.). We may also consider $[\ti X]\in\hCH_{d+1}(|W_\infty|)_\Q$ so that
we have $[W_\infty]=Z+[\ti X]$ in $\hCH_{d+1}(|W_\infty|)_\Q$. We obtain
\begin{eq}\label{trouble}
i^*(B_{|W})=B\cdot_j[W_\infty]=B\cdot_j Z+ B\cdot_j [\ti X].
\end{eq}

Denote by $\phi: |Z|\lo G$ and $\psi: \ti X\lo G$
the natural immersions.  By [GS3, Theorem 3 (4)]
the elements $B\cdot_j Z$ and $B\cdot_j[\ti X]$
are the images of the elements $B\cdot_\phi Z$ and $B\cdot_\psi [\ti X]$
of $\CH_0(|Z|)_\Q$ and $\hCH_0(\ti X)_\Q$
under the maps
$\CH_0(|Z|)_\Q\lo \hCH_0(|W_\infty|)_\Q$ and
$\hCH_0(\ti X)_\Q\lo \hCH_0(|W_\infty|)_\Q$ respectively.
We have 
$$
B\cdot_\phi Z=c_{d+1}(\xi_{||Z|})\cap [Z]
$$ 
and by [GS1, Theorem 4 (1)],
$$
B\cdot_\psi[\ti X]=\hat c_{d+1}(\bar\xi_{|\ti X})\cap [\ti X]=\hat c_{d+1}(\bar\xi_{|\ti X})
$$ 
in $\hCH_0(\ti X)_\Q=\hCH^{d+1}(\ti X)_\Q$ (recall $\ti X$ is
integral of dimension $d+1$).

Now subtract  $a(i^*_\C(\eta))$  from both sides of \ref{trouble}.
Using \ref{crap} and the above, we obtain that $i^*(b)$ can be written as a sum of
the image of the class $\hat c_{d+1}(\bar\xi_{|\ti X})-a(i^*_\C(\eta))$ 
under the map $\hCH_0(\ti X)_\Q\lo \hCH_0(|W_\infty|)_\Q$
plus the image of 
$c_{d+1}(\xi_{||Z|})\cap [Z]$ under $\CH_0(|Z|)_\Q\lo \hCH_0(|W_\infty|)_\Q$.
Since $W_\infty$ and $\ti X$ have the same generic fiber we can see from
\ref{beforecrap} that
$$
\hat c_{d+1}(\bar\xi_{|\ti X})-a(i^*_\C(\eta))=\hat c_{d+1}((\bar\zeta+A)_{|\tilde X}).
$$
Hence, part (b) will follow if we show that 
$\hat c_{d+1}((\bar\zeta+A)_{|\tilde X})=0$. 

Over $\ti X$, there is an exact sequence of vector bundles
$$
0 \lo  {\zeta_1}_{|\tilde X} \lo  {\zeta_0}_{|\tilde X}\lo Q\lo 0
$$
with $Q$ of rank $d$. We have $\tilde X_\C=X_\C$ and, as
we have seen before, 
there is an isomorphism
$Q_\C \simeq \Omega^1_{X_\C}$ which can be used to identify the above exact sequence
with $\E_\C$. This implies that 
$$
(    (\bar\zeta_0)_{|\tilde X}           , 0) -((\bar\zeta_1)_{|\tilde X}, 0)
+(0, \tch(\E_\C))=(\bar Q, 0)
$$
in $\hK_0(\ti X)$. Since $A_{|\tilde X}=A_{|X\times\{\infty\}}=(0, \tch(\E_\C))$,
this translates to 
$(\bar\zeta+A)_{|\tilde X}=(\bar Q, 0)$ in 
$\hK_0(\ti X)$. Since $d+1>{\rm rk}(Q)=d$, by [GS1, 4.9, p. 198],
$\hat c_{d+1}((\bar Q, 0))=0$. Therefore,
we obtain
$$
\hat c_{d+1}((\bar\zeta+A)_{|\tilde X})=0.
$$
This completes the proof of Lemma \ref{twores}
and therefore also of Proposition \ref{locAr}.
\medskip

\noindent{\bf Remark:} Let $\bar{\cal F}$ be a hermitian coherent
sheaf on $X$. Suppose that $Y\subset X$ is a {\sl fibral} closed subscheme 
and assume that $\F$ is locally free of rank $m$ on the complement $X-Y$. Let $z_{i,Y}: 
\CH_{d+1-i}(Y)\lo \hCH_{d+1-i}(X)$ be the natural 
homomorphism. The same argument as 
in the proof above can be used to show that for $i>m$,
$$
z_{i, Y}(c^X_{i, Y}(\F)\cap [X])=\hat c_i(\bar\F),
$$
where $c^X_{i, Y}(\F)$ is the localized
Chern class of [B, \S 1].
\medskip

\section{Tame reduction} \label{tame}

Here we show Theorem \ref{main2}. Write $I$
for an index set for the irreducible components
of the singular fibers of $X\lo \Spec(\Z)$.
If $i\in I$, we denote by $T_i$ the corresponding irreducible
component and by $m_i$ its multiplicity in
the divisor of the corresponding special fiber.
For a non-empty subset $J$ of $I$, set
$$
T_J=\cap_{i\in J}T_i
$$
(scheme-theoretic intersection). Under our assumptions,
$T_J$ is either empty or a smooth projective scheme
of dimension $d+1-|J|$ over a finite field. The union $\cup_{J\neq\subset J'} T_{J'}$
is a divisor with strict normal crossings on $T_J$.
We start with the following proposition:

\begin{prop} \label{kecoke}
With the assumptions of Theorem \ref{main2}, we can consider the
sheaf of relative logarithmic differentials $\omlog$ (see below); it is locally
free of rank $d$ on $X$. There is a morphism
$$
\omega: \Omega^1_{X/\Z}\lo \omlog
$$
whose kernel and cokernel are isomorphic
to the kernel and cokernel of the
morphism
$$
a: \oplus_{p\in S}\O_X/p\O_X\lo \oplus_{i\in S}\O_{T_{i}}.
$$
\end{prop}

\noi{\bf Proof:} The statement is local
on the base, and so to simplify notation we will
assume there is only one prime in $S$. 
We will use the logarithmic differentials $\Omega^1_{X/\Z}(\log X^{\rm red}_p)$
defined in [K] \S 2. By definition,
$$
\Omega^1_{X/\Z}(\log X^{\rm red}_p):=(\Omega^1_{X/\Z}\oplus (\O_X\otimes j_*\O^*_{X[{1\over p}]}))/\F
$$
where $j$ is the open immersion
$j: X[{1\over p}]\lo X$ and $\F$ is the $\O_X$-subsheaf generated by elements of the form
$(da, 0)-(0,a\otimes a)$ for $a\in\O_X\cap j_*\O^*_{X[{1\over p}]}$. We will 
write the element $a\otimes b$ as 
$a\cdot d\log(b)$. Notice that $j_*\O^*_{X[{1\over p}]}$ is the sheaf
of elements of the function field of $X$ whose divisor has support contained in the 
special fiber. By definition, $\omlog$ is the quotient
of $\Omega^1_{X/\Z}(\log X^{\rm red}_p)$ by the $\O_X$-subsheaf generated
by $d\log(p)$. There is an exact sequence
$$
\O_X/p\O_X\buildrel\phi\over\lo \Omega^1_{X/\Z}(\log X^{\rm red}_p)\buildrel\omega_1\over\lo \omlog\lo 0
$$
where the homomorphism $\phi$
maps $f$ to $f\cdot d\log (p)$.
There is also a natural exact sequence
\begin{eq}\label{seque}
0\lo \Omega^1_{X/\Z}\buildrel\omega_2\over\lo \Omega^1_{X/\Z}(\log X^{\rm red}_p)\buildrel \oplus_i {\rm Res_i}\over \lo \oplus_i\O_{T_i}\lo 0.
\end{eq}
Here the right hand homomorphism is given by taking residues
along $T_i$. The homomorphism $\omega$ is equal to
the composition $\omega_1\cdot\omega_2$.

Under our assumptions, the scheme $X$ is locally \'etale
isomorphic to
$$
Y=\Spec(\Z[t_1, \ldots ,t_{d}]/(t_1^{m_1}\cdots t_d^{m_d}-p))
$$
with all $m_i$ prime to $p$. The above constructions
of logarithmic differentials etc. make sense for the scheme $Y$;
we can see by an explicit calculation that
$\phi_Y$ is injective and that the analogue of the 
sequence \ref{seque} for $Y$ is exact. It follows from the fact that
taking (logarithmic) differentials commutes with \'etale base change
that $\phi$ is injective and that the 
sequence \ref{seque} is exact.
On $Y$ we have $t^{m_1}_1\cdots t_d^{m_d}=p$ and so 
$$
d\log(p)=m_1{dt_1\over t_1}+\cdots +m_d{dt_d\over t_d}.
$$
This shows that for $f\in \O_Y/p\O_Y$, $\phi_Y(f)$ gives an element in the kernel
of $\omega$ if and only if $f\in (t_1\cdots t_d)$; 
this translates to $a(f)=0$. Furthermore, $f\cdot d\log(p)=0$ if and only if $f=0$
in $\O_Y/p\O_Y$. This shows 
the statement about the kernels
for $X$. 
Let us now discuss the cokernels: Let $\beta:\oplus_i\O_{T_i}\lo \oplus_i\O_{T_i}$
be the automorphism defined by $\beta((f_i)_i)=(m_if_i)_i$
(recall that all the $m_i$ are prime to $p$).
The above calculation on $Y$ implies that the composition
$$
\O_X/p\O_X\buildrel\phi \over\lo \Omega^1_{X/\Z}(\log X^{\rm red}_p) 
\buildrel \oplus_i {\rm Res_i}\over \lo \oplus_i\O_{T_i}
$$
coincides with $f\mapsto (m_1f,\cdots, m_df)$. The residue homomorphism ${\rm Res}=\oplus_i{\rm Res}_i$
now gives
a surjection:
$$
\omlog\buildrel\beta^{-1}\cdot {\rm Res}\over\lo {\rm coker}(a)\lo 0
$$
and we have ${\rm ker}(\beta^{-1}\cdot {\rm Res})={\rm ker}({\rm Res})=\omega(\Omega^1_{X/\Z})$.
This implies
${\rm coker}(\omega)\simeq {\rm coker}(a)$.
\end{Proof}
\smallskip
\medskip

Let $K_0^{X_S}(X)$ be the Grothendieck group
of complexes of
locally free $\O_X$-sheaves which are exact
off $X_S$; since $X$ is regular, $K_0^{X_S}(X)$
can be identified with $K_0'(X_S)$.
Set $q=\prod_{p\in S}p$.
Consider the following complexes of
locally free $\O_X$-sheaves which are exact
off $X_S$:
$$
\E_1: N\buildrel \delta\over \lo P\lo \omlog
$$
$$
\E_2: \O_X\buildrel  (q, -q)\over\lo \O_X\oplus(\oplus_i \O_X(-T_i))\lo \oplus_i\O_X
$$
concentrated in degrees $-1,0,1$. The second
homomorphism of $\E_1$ is the composition 
of $P\lo \Omega^1_{X/\Z}$ with $\omega$; the 
second homomorphism of $\E_2$ is given by $(g, (h_i)_i)\mapsto (g+h_i)_i$.  
Proposition \ref{kecoke} implies that $[\E_1]=[\E_2]$ 
in $K_0^{X_S}(X)$. 
Consider also the complex
$$
\E_3: N\buildrel (\delta,0)\over \lo P\oplus\omlog\buildrel (0, id)\over\lo \omlog
$$
concentrated in degrees $-1,0,1$. The complex $\E_3$ is quasi-isomorphic to 
the complex $N\buildrel \delta\over\lo P$ (in degrees $-1$ and $0$).
There is an exact sequence of complexes
$$
0\lo \E_1\lo \E_3\buildrel pr\over \lo \omlog\lo 0
$$
where on the right end, $\omlog$ is considered as a complex supported on 
degree $0$. Therefore, the main result of [A] (see loc. cit. Proposition 1.4 also [B] Prop. 1.1)
implies that
\begin{eq}\label{fult}
c_{d+1}^{X_S}(\Omega^1_{X/\Z})=\sum_{k+l=d+1}c_{k}(\omlog)\cdot c^{X_S}_l([\E_1]).
\end{eq}
In fact, since $[\E_1]=[\E_2]$ we can replace $c^{X_S}_l([\E_1])$
by $c^{X_S}_l([\E_2])$ in this equality.
We have
$$
[\E_2]=[\O_X/q\O_X]-\sum_i[\O_{T_i}]
$$
(here we identify $K_0^{X_S}(X)$ with $K'_0(X_S)$)
and so
\begin{eq}\label{ce1}
c_l^{X_S}([\E_2])=c_l^{X_S}([\O_X/q\O_X]+\sum_i(-[\O_{T_i}])).
\end{eq}
We have $c^{X_S}_1([\O_X/q\O_X])=\sum_i m_i[T_i]$, $c^{X_S}_l([\O_X/q\O_X])=0$
for $l>1$. Similarly, $c^{X_S}_1(-[\O_{T_i}])=-[T_i]$, $c^{X_S}_l(-[\O_{T_i}])=0$, for $l>1$.
Combining these with \ref{ce1} we obtain from the usual Chern class identities
\begin{eq} \label{ce2}
c_l^{X_S}([\E_2])=\sum_{J\subset I, |J|=l}(-1)^{|J|}[T_J]+
(\sum_{i\in I}m_i[T_i])(\sum_{J'\subset I, |J'|=l-1}(-1)^{|J'|}[T_{J'}]).
\end{eq}
Now since $\sum_im_i[T_i]$ is a principal divisor
in $X$ we get for $l\geq 2$ 
$$
(\sum_{i\in I}m_i[T_i])(\sum_{J'\subset I, |J'|=l-1}(-1)^{|J'|}[T_{J'}])=0\in \CH_*(X_S).
$$

Combining this with \ref{fult} and \ref{ce2} we get
\begin{eqnar} \label{chern}
c_{d+1}^{X_S}(\Omega^1_{X/\Z})&=&\sum_{i\in I}(m_i-1)c_d(\omlog)\cdot [T_i]+\\ \nonumber
&&+\sum_{J\subset I, |J|\geq 2}(-1)^{|J|}c_{d+1-|J|}(\omlog)\cdot [T_J].
\end{eqnar}
Therefore
\begin{eqnar}\label{chern2}
c_{d+1}^{X_S}(\Omega^1_{X/\Z})&=&\sum_{i\in I}(m_i-1)c_d(\omlog_{|T_i})+\\ \nonumber
&&+\sum_{J\subset I, |J|\geq 2}(-1)^{|J|}c_{d+1-|J|}(\omlog_{|T_J}).
\end{eqnar}

\begin{prop} \label{euler}
For a non-empty subset $J$ of $I$, set $T^*_J=T_J-\cup_{J\neq\subset J'} T_{J'}$.
We have
$$
{\rm deg}(c_{d+1-|J|}(\omlog_{|T_J}))=(-1)^{d+1-|J|}\chi_c(T^*_J)
$$
where $\chi_c(T^*_J)$ is the $l$-adic ($l\notin S$)
Euler characteristic with compact supports of $T^*_J$.
\end{prop}

\begin{Proof}
Denote by $\log X^{\rm red}_p|T_J$
the logarithmic structure on $T_J$ obtained by restricting the logarithmic structure
given by $(X, X^{\rm red}_p)$ to $T_J$. This is isomorphic to the
logarithmic structure defined on $T_J$ by its divisor
with strict normal crossings $\cup_{J\neq\subset J'} T_{J'}$. We will show
that
\begin{eq}\label{show}
[\omlog_{|T_J}]=[\Omega^1_{T_J/k}(\log X^{\rm red}_p|T_J)]+(|J|-1)[\O_{T_J}]
\end{eq}
in ${\rm K}_0(T_J)$. The proposition will follow
from \ref{show} and the well-known fact (see for example
[S], p. 402) that 
$$
\deg (c_{d+1-|J|}(\Omega^1_{T_J/k}(\log X^{\rm red}_p|T_J)))=(-1)^{d+1-|J|}\chi_c(T^*_J).
$$
From the proof of Proposition \ref{kecoke}
there is an exact sequence
\begin{eq}\label{exact1}
0\lo \O_{T_i}\lo \Omega^1_{X/\Z}(\log X^{\rm red}_p)_{|T_i}\lo \omlog _{|T_i}\lo 0.
\end{eq}
By [K] \S2 (see also [S], p. 404) there are also
exact sequences
\begin{eq}\label{exact2}
0\lo \Omega^1_{T_i/{\bf F}_p}(\log X^{\rm red}_p|T_i )\lo
\Omega^1_{X/\Z}(\log X^{\rm red}_p)_{|T_i}\lo \O_{T_i}\lo 0,
\end{eq}
and for $|J'|=|J|+1$,
\begin{eq}\label{exact2}
0\lo \Omega^1_{T_{J'}/{\bf F}_p}(\log X^{\rm red}_p|T_{J'} )\lo
\Omega^1_{T_J/{\bf F}_p}(\log X^{\rm red}_p|T_{J})_{|T_{J'}}\lo \O_{T_{J'}}\lo 0.
\end{eq}
We can now see that \ref{show} follows by induction on the cardinality of $J$.

\end{Proof}

\medskip
\medskip

Proposition \ref{euler} and \ref{chern2} give for $p\in S$:

\begin{eqnarray} \label{equ}
\deg((-1)^{d+1}c_{d+1}^{X_S}({\Omega^1_{X/\Z}})_{|X_p})   &=&-\sum_{i\in {I_p}}(m_i-1)
\chi^*_c(T_i)+\sum_{J\subset I_p, |J|\geq 2}\chi^*_c(T_J)\nonumber\\
&=&-\sum_{i\in {I_p}} m_i\chi^*_c (T_i)+\chi(X_p)
\end{eqnarray}
where $I_p$ is the subset of $I$ that
corresponds to components over $p$.

Under our assumption, the ramification is tame 
(there is no Swan term in the conductor) and
for each $p\in S$,
$$
\chi(X_\Q)=\sum_{i\in I_p}m_i\chi^*_c(T_i)
$$
(see for example [S], Cor. 2, p. 407). Therefore, 
$$
A(X)=\prod_{p\in S} p^{\chi(X_\Q)-\chi(X_p)}=\prod_{p\in S} p^{\sum_i m_i\chi^*_c(T_i)-\chi(X_p)}.
$$
This together with \ref{equ} completes the proof of \ref{main2}.

\bigskip

\noindent{\sc T. Chinburg, \  University of Pennsylvania, Phila., PA
19104.}

 ted@math.upenn.edu

\medskip

\noindent{\sc G. Pappas, \ Michigan State University, E. Lansing, MI 48824.} 

 pappas@math.msu.edu
\medskip

\noindent{\sc M. J. Taylor, \ UMIST, Manchester, M60 1QD, UK.}

 Martin.Taylor@umist.ac.uk

\end{document}